\documentstyle{article}

\newtheorem{th}{Theorem}
\newtheorem{lem}{Lemma}

\begin{document}
\author{A.A.Davydov}
\title{On some Hochschild cohomology classes of fusion algebras}
\maketitle
\begin{abstract}
The obstructions for an arbitrary fusion algebra to be a fusion algebra of some 
semisimple monoidal category are constructed. Those obstructions lie in groups 
which are closely related to the 
Hochschild cohomology of fusion algebras with coefficients in the $K$-theory 
of the ground (algebraically closed) field.

The special attention is devoted to the case of fusion algebra of invariants 
of finite group action on the group ring of abelian group. 
\end{abstract}
\tableofcontents
\section{Introduction}
Recent activity in the theory of fusion rings (a special kind of semiring) 
was inspired by the development of quantum field theory \cite{ver,ms,wit}, 
although this notion appeared before in different branches of mathematics, 
for example in the representation theory (semiring of irreducible 
representations). Most fusion rings considered so far can be 
represented as fusion rings of simple objects of some semisimple monoidal 
category \cite{dm} (so-called monoidal fusion rings). As we shall see not all 
fusion rings admit such representation. 

The article is devoted to the construction of obstructions for an arbitrary 
fusion ring to be the monoidal. These obstrucions are elements of the Hochschild 
cohomology of the fusion ring with coefficients in the algebraic K-theory of the 
ground field. The construction is based on the notion of $A_n$-space (homotopy 
associative space) \cite{sta,may1,may2,bv}. The obstruction to monoidality is a 
special case of the obstruction to extension of $A_n$-structure to $A_{n+1}$. 

The first obstruction to monoidality, which is an element of the fourth 
Hochschild cohomology can be expressed explicitly by means of structural 
constants of the given fusion ring (section \ref{com}). In section \ref{calc} 
the analysis of triviality of the first obstruction for the case of two element 
fusion ring with identity. 

\section*{Acknowledgement}

The author would like to thank J.Stasheff for stimulating attention, 
fruitful discussions and corrections.

\section{Fusion rings}
A {\em fusion ring} is a set $S$ with a collection of non-negative integers 
$\{ m^{x}_{x_1 ,x_2}, x,x_1 ,x_2\in S\}$ ({\em structural constants}) which 
satisfy the ({\em associativity}) condititon
$$m^{x}_{x_1 ,x_2 ,x_3} = \sum_{t\in S}m^{x}_{x_1 ,t}m^{t}_{x_2 ,x_3} = 
\sum_{s\in S}m^{s}_{x_1 ,x_2}m^{x}_{s,x_3}, \quad \forall x,x_1 ,x_2 ,
x_3\in S.$$
An element $e$ of the fusion algebra $S$ is a {\em identity} if $m^{s}_{t,e} = 
m^{s}_{e,t} = \delta_{s,t}$ for all $s,t\in S$. 

A {\em morphism} of the fusion ring $S$ to the fusion ring $S'$ is a 
collection $\{ n^{s}_{t}, s\in S, t\in S' \}$ of non-negative integers which 
satisfy the following condition:
$$\sum_{s\in S}m^{s}_{s_1 ,s_2}n^{t}_{s} = 
\sum_{t_1 ,t_2 \in S'}{m'}^{t}_{t_1 ,t_2}n^{t_1}_{s_1}n^{t_2}_{s_2}$$
for any $s_1 ,s_2 \in S$ and $t\in S'$. 

The {\em enveloping ring} $A(S)$ of the fusion algebra $S$ is the free 
$\bf Z$-module with the basis $\{[s], s\in S\}$ labeled by the elements of 
$S$ and with the multiplications $[i][j] = \sum_{s\in S}m^{s}_{i,j}[s]$. 
A morphism of fusion rings defines a homomorphism of their enveloping 
rings $f:A(S)\to A(S')$ where $f([s]) = \sum_{t\in S'}n^{t}_{s}[t]$. 

\section{Semisimple monoidal categories}
This section is devoted to the most important case of fusion algebras, namely 
fusion algebras corresponding to semisimple monoidal categories 
(monoidal fusion algebras). 

Let $\cal G$ be a semisimple monoidal $k$-linear category over the field $k$ 
with the set $S$ of isomorphism classes of simple objects. The collection of
dimensions $m^{x}_{y,z} = dim_k Hom_{\cal G}(X,Y\otimes Z)$ form a fusion ring 
structure on the set $S$. Here $X,Y$ and $Z$ are some representatives of the 
classes $x,y,z\in S$. This fusion ring structure (together with the map from
$S$ to the set of simple $k$-algebras, $x\mapsto End_{\cal G}(X)$) contains all 
the information about the tensor product in the category $\cal G$. But as we 
shall see it is not sufficient to reconstruct associativity constraint $\phi$ of 
the category. To define this constraint in terms of simple objects it is 
necessary to fix representative of any class from 
$S$. If we do this we can construct a semisimple monoidal category in which 
isomorphic simple objects coincides and which is equivalent to the category 
$\cal G$. We will call such category a {\em model} ${\cal M}(S)$ of the 
semisimple monoidal category $\cal G$ with the set $S$ of simple objects. 
Objects of the model are maps from the set $S$ to the set of finite dimensional 
vector spaces with only finitely many nonzero values (or formal finite sums 
$\oplus_{x\in S}V_{x}x$ of vector spaces).
The morphism space $Hom(\oplus_{x\in S}V_{x}x,\oplus_{y\in S}U_{y}y)$ coincides 
with $\oplus_{x\in S}Hom_k (V_{x},U_{x})\otimes End_{\cal G}(X)$. The tensor 
product of two simple objects $x\otimes y$ is by definition the sum 
$\oplus_{z\in S}H^{z}_{x,y}$, where $H^{z}_{x,y} = Hom_{\cal G}(Z,X\otimes Y)$.  
The associativity constraint for the model can be given by the collection of 
isomorphisms of vector spaces 
\begin{equation}\label{asc}
\Phi^{x}_{x_1 ,x_2 ,x_3}:\oplus_{u\in S}H^{x}_{x_1 ,u}\otimes H^{u}_{x_2 ,x_3}
\to\oplus_{v\in S}H^{v}_{x_1 ,x_2}\otimes H^{x}_{v,x_3},
\end{equation}  
which are defined by the associativity constraint $\phi$ of the category 
$\cal G$ 
$$\Phi^{x}_{x_1 ,x_2 ,x_3} = 
Hom_{\cal G}(X,\phi_{X_1 ,X_2 ,X_3}):Hom_{\cal G}(X,X_1\otimes (X_2\otimes X_3 
))\to Hom_{\cal G}(X,(X_1\otimes X_2 )\otimes X_3 ).$$
The pentagon axiom for the constraint $\phi$ consists of commutativity of the 
pentagon diagram \cite{dm}
for any $x,x_2 ,x_3 ,x_3 ,x_4\in S$. This condition can be written in the form 
of an equation: 
$$(\oplus_{b\in S}(\Phi^{x}_{b,x_3 ,x_4})_{23})
(\oplus_{a\in S}(\Phi^{t}_{x_1 ,x_2 ,a})_{12}) = $$
$$(\oplus_{f\in S}(\Phi^{f}_{x_1 ,x_2 ,x_3})_{12})
(\oplus_{d\in S}(\Phi^{t}_{x_1 ,d,x_4})_{13})
(\oplus_{c\in S}(\Phi^{c}_{x_2 ,x_3 ,x_4})_{23}),$$  
which will be denoted $A^{x}_{x_1 ,x_2 ,x_3 ,x_4}$. 

Later for the case of simplicity we will consider $k$-linear semisimple 
categories for which $End(X) = k$ for any simple object $X$. For example, it is 
true for algebraically closed field $k$. 

If we start from arbitrary fusion ring $S$ we also can construct the semisimple 
category ${\cal G}(S)$ with tensor product whose simple objects are parametrized 
by the elements of $S$. Since the vector spaces 
\begin{equation}\label{hs}
\oplus_{u\in S}H^{x}_{x_1 ,u}\otimes H^{u}_{x_2 ,x_3}\quad\mbox{and}\quad
\oplus_{v\in S}H^{v}_{x_1 ,x_2}\otimes H^{x}_{v,x_3}
\end{equation}
have the same dimension $m^{x}_{x_1 ,x_2 ,x_3}$, they are isomorphic, i.e. the 
tensor product is {\em quasiassociative}. It is associative if we can choose the 
isomorphisms between vector spaces (\ref{hs}) satisfying the pentagon axiom. 

In the case of the fusion ring consisting of one element $x$ with the fusion 
rule $x*x = nx$ (the structural constant $m^{x}_{x,x}$ equals $n$) the 
quasiassociative structure is an automorphism $\Phi$ of vector space 
$H\otimes H$. Here $H = H^{x}_{x,x}$ is an $n$-dimensional vector space. The 
pentagon axiom for this structure is equivalent to the following 
({\em pentagon}) equation
\begin{equation}\label{pe}
\Phi_{12}\Phi_{13}\Phi_{23} = \Phi_{23}\Phi_{12}.
\end{equation}  
Here $\Phi_{ij}$ denotes the automorphism of $H^{\otimes 3}$ acting on i-th and 
j-th components. 

To any structure of Hopf algebra on the vector space $H$, we can associate a 
solution to the equation (\ref{pe}):  
$$\Phi (g\otimes h) = \sum_{(g)}g_{(0)}\otimes h_{(0)}\quad\forall g,h\in H.$$
General solutions to equation (\ref{pe}) are very close to these. The 
proof of this fact will be published elsewhere. 

The pair consisting of the identity $e$ and the element $x$ with fusion 
rule $x*x = ne$ provides an example of a fusion ring whose category does not 
admit a monoidal structure (there is no quasiassociative structure satisfying 
the pentagon axiom). To simplify the verification of this fact let us prove the 
following lemma.
\begin{lem}
Let $S$ be a fusion ring with identity $e$ and $\Phi$ be an associativity 
constraint for a corresponding semisimple $k$-linear category. Then there are 
an identifications $\rho_{x}:H^{x}_{x,e}\to k$ and 
$\lambda_{x}:H^{x}_{e,x}\to k$ (for any $x$) such that $\Phi^{z}_{x,y,e}, 
\Phi^{z}_{x,e,y}$ and $\Phi^{z}_{e,x,y}$ coincide with the compositions: 
\begin{equation}\label{fe}
(I\otimes\rho_{z}^{-1})(I\otimes\rho_{y}):H^{z}_{x,y}\otimes H^{y}_{y,e}\to 
H^{z}_{x,y}\to H^{z}_{x,y}\otimes H^{z}_{z,e},
\end{equation}
\begin{equation}\label{se}
(\rho_{x}^{-1}\otimes I)(I\otimes\lambda_{y}):H^{z}_{x,y}\otimes 
H^{y}_{e,y}\to H^{z}_{x,y}\to H^{x}_{x,e}\otimes H^{z}_{z,e}, 
\end{equation}
\begin{equation}\label{te}
(\lambda_{x}^{-1}\otimes I)(\lambda_{z}\otimes I):H^{z}_{e,z}\otimes 
H^{z}_{x,y}\to H^{z}_{x,y}\to H^{x}_{e,x}\otimes H^{z}_{z,e}.
\end{equation}
\end{lem}
Proof: 

It follows from the equations $A^{x}_{x,e,e,e}, A^{x}_{e,x,e,e}, 
A^{x}_{e,e,x,e}$ and $A^{x}_{e,e,e,x}$ that we can choose the isomorphisms 
$\rho_{x}:H^{x}_{x,e}\to k$ and $\lambda_{x}:H^{x}_{e,x}\to k$ such that 
$\Phi^{x}_{x,e,e}, \Phi^{x}_{e,x,e}$ and $\Phi^{x}_{e,e,x}$ will have the form 
\ref{se},\ref{fe} and \ref{te} respectively. 

Conversely equations which involve three identity elements imply the statement 
of the lemma. 
$\Box$

Let us return to the fusion ring $\{ e,x\}$ with fusion rule $x*x = ne$. 
Identifying $H^{x}_{x,e}$ and $H^{x}_{e,x}$ with $k$ we can see that unique 
nontrivial component $\phi^{x}_{x,x,x}:H^{x}_{x,e}\otimes H^{e}_{x,x}\to 
H^{e}_{x,x}\otimes H^{x}_{e,x}$ of a monoidal structure for the category 
${\cal G}(S)$ (if it exists) corresponds to automorphism $\Phi$ of 
$n$-dimensional vector space $H=H^{e}_{x,x}$. The pentagon condition 
$A^{e}_{x,x,x,x}$ is equivalent to the equation $\Phi^2 \otimes I = t$, 
where $I$ is the identity automorphism of the space $H$ and $t:H^{\otimes 2}\to 
H^{\otimes 2}$ is the permutation of tensor factors (the isomorphisms 
$(\Phi^{e}_{x,x,e})_{12}$ and $(\Phi^{e}_{e,x,x})_{23}$ 
are identical by the lemma, $(\Phi^{x}_{x,x,x})_{23}$ and 
$(\Phi^{x}_{x,x,x})_{12}$ correspond to $I\otimes A$ and $A\otimes I$ 
respectively and $(\Phi^{e}_{x,e,x})_{13}$ to $t$). 
It is easy to see that this equation has solutions iff $n=1$.

This example motivates the following definition. A fusion ring will be 
called {\em monoidal} if it coincides with the fusion ring of simple objects of 
some semisimple monoidal category. 

Any monoidal $k$-linear functor \cite{dm} $F:{\cal G}\to{\cal G}'$ between 
semisimple 
monoidal $k$-linear categories provides the morphism of fusion algebras $S$ and 
$S'$ of simple objects. This morphism is given by the collection of dimensions 
$n^{x}_{y} = dim_k Hom_{{\cal G}'}(X,F(Y))$. Here $X\in{\cal G}'$ and 
$Y\in\cal G$ are some representatives of the classes $x\in S'$,$y\in S$. This 
morphism structure contains all information about the functor $F$ (without its 
monoidal structure). 
Let us denote by $H^{x}_{y}$ the $n^{x}_{y}$-dimensional vector space 
$Hom_{{\cal G}'}(X,F(Y))$. 

The monoidal constraint for this functor can be given by the collection of 
isomorphisms of vector spaces 
$$\Psi^{x}_{y,z}:\oplus_{u\in S}H^{u}_{y,z}\otimes H^{x}_{u}\to
\oplus_{v,w\in S'}H^{v}_{y}\otimes H^{w}_{z}\otimes H^{x}_{v,w},$$  
which are defined by the monoidal constraint $\psi$ of the functor $F$: 
$$\Psi^{x}_{y,z} = 
Hom_{{\cal G}'}(X,\phi_{Y,Z}):Hom_{{\cal G}'}(X,F(Y\otimes Z))\to 
Hom_{{\cal G}'}(X,F(Y)\otimes F(Z)).$$
The compatibility axiom for the constraint $\phi$ consists of commutativity 
of hexagon diagrams \cite{dm}
for any $x\in S',y,z,w\in S$. This condition can be written in the form of an 
equation: 
$$(\oplus_{t_1 ,t_2 ,t_3 \in S'}(\Phi^{t}_{t_1 ,t_2 ,t_3})_{13})
(\oplus_{y\in S'}(\Psi^{y}_{s_1 ,s_2})_{12})
(\oplus_{x\in S}t_{12}(\Psi^{t}_{s,x})_{23})t_{12}) = $$
$$t_{34}(\oplus_{w\in S'}(\Psi^{w}_{s_1 ,s_2})_{12})
(\oplus_{z\in S}(\Psi^{t}_{z,s_3})_{23})
(\oplus_{s\in S}(\Phi^{s}_{s_1 ,s_2 ,s_3})_{12}).$$  
The condition of existence of a monoidal functor between semisimple monoidal 
categories ${\cal M}(S)$, ${\cal M}(S')$ which induce a given morphism of fusion 
rings $S,S'$ will be called {\em monoidality} of the morphism (with respect to 
the given monoidal structures on the categories ${\cal M}(S)$, ${\cal M}(S')$). 

The simplest example of a nonmonoidal morphism of fusion rings is provided by 
the identity map of the one-element fusion ring $S = \{ x, x*x=nx\}$ with 
respect to monoidal structures on the category ${\cal M}(S)$ defined by two 
non-isomorphic $n$-dimensional Hopf algebras.

\section{$A_n$-spaces}

Let us denote by $T^k$ the set of planar trees with $k$ ends. 
This set is a disjoint union $\bigsqcup_{i=0,...,k-2}T^{k}_{i}$ of subsets 
$T^{k}_{i}$, which consist of trees with $k-i-1$ vertices. Junction of two 
neighboring vertices 
defines the 
correspondence $R^{k}\subset T^{k}\times T^{k}$. Namely, 
$$R^k =\{ (t,t')\in T^{k}\times T^{k}, \mbox{where} t' 
\mbox{is a result of the contraction of an internal edge of} t\}.$$
This correspondence is also 
decomposes into a disjoint union $\bigsqcup_{i=0,...,k-3}R^{k}_{i}$, where 
$R^{k}_{i}$ is contained in $T^{k}_{i}\times T^{k}_{i+1}$. 
The unique element of $T^{k}_{k-2}$ will be denoted by $t(k)$. 

This data allows us to define cell complexes $BT^{k}$ (Stasheff's complexes 
$K_{k}$,
\cite{sta}): elements of the set $T^{k}_{i}$ parametrizes the $i$-dimensional 
cells of $BT^{k}$, which are glued by means of correspondence $R$.

The operation of gluing trees (roots of the trees $t_1 ,...,t_k$ glued together 
with the ends of the tree $t\in T^k$) defines the maps 
$$c^{k,l_1 ,...,l_k}:T^{k}_{i}\times T^{l_1}_{i_1}\times ...\times T^{l_k}_{i_k}
\to T^{l_1 +...+l_k}_{i_1 +...+l_k +i}.$$
These maps induce a (non-symmetric) operad structure \cite{may1}:
$$BT^{k}\times BT^{l_1}\times ...\times BT^{l_k}\to BT^{l_1 +...+l_k}.$$ 
Let us denote by $\partial^{i}_{k}(t)\in T^{k+1}$ the result of gluing the root 
of the tree $t(2)$ to i-th end of the tree $t$ for $t\in T^k$ and $i=1,...,k$. 
For $i=0,k+1$ the expression $\partial^{i}_{k}(t)$ will denote the 
result of gluing the root of the tree $t$ to the second and first ends of the 
tree $t(2)$ respectively. 

Let us shoose the orientations of the complexes $BT^k$ such that the maps
$\partial^{1}_{k}:BT^{k}\to BT^{k+1}$ will preserve the orientations. 

It should be noted that for any space $X$ the collection of spaces $Hom(X^k ,X)$ 
also forms an operad. 

An {\em $A_n$-space} structure \cite{sta,bv,may1} is a collection of continuous 
maps $BT^k\to Hom(X^k ,X)$ for $k\leq n$, which are compatible with 
compositions. 

It is easy to see that the cells which correspond to the trees from the 
image of the composition map 
$$c^{k,l_1 ,...,l_k}(T^{k}_{i}\times T^{l_1}_{i_1}\times 
...\times T^{l_k}_{i_k})$$
are contractible in $Hom(X^{l_1 +...+l_k},X)$ if $i_p , i_q >0$ for some $p,q$. 

For example an $A_3$-space structure on $X$ consists of a product $\mu 
:X\times X\to X$ (the image of $pt=BT^2 \to Hom(X^2 ,X)$) and a homotopy 
between $\mu (I\wedge\mu)$ and $\mu (\mu\wedge I)$ (the image of $I=BT^3 \to 
Hom(X^3 ,X)$). 

A $A_{n+i}$-structure on the space $X$ for which the restriction on the 
$k$-cells for $k\leq n$ coincides with the given $A_n$-structure will be 
called an {\em extension}. Two $A_n$-extensions of a given $A_{n-1}$-structure 
on the space $X$ define a map 
$X^{n}\to\Omega^{n-2}X$. Indeed these $A_n$-structures differ only by the maps
$I^{n-2} = BT^{n}\to Hom(X^{n},X)$ of the $n-2$-cell. In particular these maps 
coincide on the boundary $S^{n-3}=I^{n-2}$. Hence we can glue them together 
along the boundary and create a map $S^{n-2} = 
I^{n-2}\bigsqcup_{\partial I^{n-2}}I^{n-2}\to Hom(X^{n},X)$ or 
$X^{n}\to\Omega^{n-2}X$. 

Now let us define the cohomology of $A_3$-spaces with coefficients in some 
bimodule.

An {\em $A_3$-bimodule} over an $A_3$-space $X$ is a space $Y$ 
together with the continuous maps 
$$\nu :Y\wedge X\rightarrow Y, \qquad\upsilon :X\wedge Y\rightarrow Y$$
and homotopies
$$\nu (I\wedge\mu )\rightarrow\nu (\nu\wedge I), \quad
\upsilon (I\wedge\nu )\rightarrow\nu (\upsilon\wedge I), \quad \upsilon 
(I\wedge\upsilon )\rightarrow\upsilon (\mu\wedge I).$$
For example the space of maps to the $A_3$-space is a $A_3$-bimodule over this 
space.

Denote by $[X,Y]$ the set of homotopy classes of continuous maps from $X$ to 
$Y$.
The {\em topological cobar} complex of a $A_3$-space $X$ with coefficients in a 
$A_3$-bimodule space $Y$ is a cosimplicial complex of sets
$$C_{*}(X,Y), \qquad C_{n}(X,Y) = [X^{n},Y]$$
with the coface maps ${\partial}^{i}_{n}:C_{n-1}(X,Y)\rightarrow C_{n}(X,Y)$ 
defined as follows
$${\partial}^{i}_{n}(f) = \left\{ \begin{array}{l}
\upsilon (I\times f),\qquad i=0\\
f(I\times ...\times\mu\times ...\times I),\qquad 1\leq i\leq n\\
\nu(f\times I),\qquad i=n+1
\end{array}\right.$$
If $Y$ is a loop space (with a loop structure which is not connected with the 
given $A_3$-structure) then $C_{\cdot}(X,Y)$ is a complex of groups and these 
groups are abelian if $Y$ is a double loop space. In the second case the 
cohomology of the cochain complex associated with $C_{*}(X,Y)$ will be called 
the {\em topological} cohomology ($H^{*}(X,Y)$) of the $A_3$-space $X$ with 
coefficient in the $A_3$-bimodule $Y$.

Let us note that an $A_n$-space structure on $X$ (and the operad structures on 
$BT^*$ and $Hom(X^* ,X)$) allows to define maps 
$BT^{l}_{\leq n-2}\to Hom(X^l ,X)$ for any $l$. Here $BT^{l}_{\leq n-2}$ denotes 
the union of $i$-dimensional cells in $BT^l$ with $i\leq n-2$. 

For example the map $S^{n-2}=\partial BT^{n+1}=BT^{n+1}_{\leq n-1}\to 
Hom(X^{n+1},X)$ defines an element in $Hom(X^{n+1},\Omega^{n-2}X)$. 
\begin{th}\label{space}
Let $X$ be an $A_n$-space for $n\geq 3$ (and independently a loop space if 
$n=3$). Then the class in 
$[X^{n+1},\Omega^{n-2}X]$ of the map $BT^{n+1}_{\leq n-2}\to Hom(X^{n+1},X)$ is 
a cocycle. Its class $\alpha (X)$ in cohomology $H^{n+1}(X,\Omega^{n-2}X)$ does 
not depend of the choice of the $n-2$-dimesional component of $A_n$-structure on 
the space  $X$ and is trivial iff there is a $n-2$-deformation of the given 
$A_n$-structure which can be extended to an $A_{n+1}$-structure.
\end{th}
Proof:

Let us consider the map $BT^{n+2}_{\leq n-2}\to Hom(X^{n+2},X)$. The complex  
$BT^{n+2}_{\leq n-2}$ is a union of $n-2$-spheres (boundaries of $n-1$-cells), 
which correspond to planar trees with $n+2$ ends and two vertices. Hence this 
map provides an equation on the classes of these spheres in 
$[X^{n+2},\Omega^{n-2}X]$.  
The cells labeled by trees, whose vertices have valences (the number of incident, 
e,g, incoming and outcoming edges) greater then three, are 
contractible in $Hom(X^{n+2},X)$. So (potentially) non-trivial components 
correspond to the trees which can be presented as glueing of trees $t(n+1)$ and 
$t(2)$, e.g trees of the form $\partial^{i}_{n+1}(t)$ for $i = 0,...,n+2$. 
Their classes in $[X^{n+2},\Omega^{n-1}X]$ are equal to 
$\partial^{i}_{n+1}(\alpha (X))$ respectively. Hence the map 
$BT^{n+2}_{\leq n-2}\to Hom(X^{n+2},X)$ defines the equation $\partial 
(\alpha (X)) = 0$ in $[X^{n+2},\Omega^{n-2}X]$ which means that $\alpha (X)$ is a 
$n+1$-cocycle. 

Another choice of the $n$-component of an $A_n$-structure provides the cocycle 
which differs from the given one by coboundary $\partial (f)$ of the map $f$ 
which is defined by two $A_n$-extensions of the given $A_{n-1}$-structure. 

Triviality of the cohomology class of the cocycle $\alpha$ means that we can 
choose the $n$-component of $A_n$-structure such that it can be extendable to 
an $A_{n+1}$-structure. 
$\Box$

To define $A_n$-maps between $A_n$-spaces we need the following notion. 
A {\em tricolored} tree is a planar tree which vertexes are labeled by one of 
three ordered indices $0<\epsilon <1$. The index of vertix which lies over are 
not less then index of given vertix, moreover the index of the vertix which lies 
over $\epsilon$-labeled vertix is $1$. 

As in the non-colored case the set $CT^k$ of tricolored trees with $k$ 
ends can be presented as the disjoint union $\bigsqcup_{i=0,...,k-1}CT^{k}_{i}$ 
of subsets, which consist of trees with $k-i-1+n_{\epsilon}$ vertices. Here 
$n_{\epsilon}$ is the number of $\epsilon$-labeled vertices. 

We also define the color version $CR$ of the correspondence $R$ by the condition 
that $(t,t')\in CR$ if the tree $t'$ can be obtained from the tree $t$ by 
contraction of an edge between equally labeled vertices or between vertices one 
of which is labeled by $\epsilon$. 

This allows us to define cell complexes $BCT^{k}$ whose $i$-cells are 
parametrized by elements of $CT^{k}_{i}$ and are glued by means of the 
correspondence $CR$.

Two collection of inclusions $T^{k}_{i}\to CT^{k}_{i}$ which label all vertices 
by $0$ and $1$ respectively provides two collections of compositions:
$$u^{k,l_1 ,...,l_k}:T^{k}_{i}\times CT^{l_1}_{i_1}\times ...\times 
CT^{l_k}_{i_k}\to CT^{l_1 +...+l_k}_{i_1 +...+l_k +i},$$
and
$$d^{k,l_1 ,...,l_k}:CT^{k}_{i}\times T^{l_1}_{i_1}\times ...\times 
T^{l_k}_{i_k}\to CT^{l_1 +...+l_k}_{i_1 +...+l_k +i},$$
which satisfy the usual operad axioms and define two collections of cell 
complex maps:
$$u^{k,l_1 ,...,l_k}:BT^{k}\times BCT^{l_1}\times ...\times BCT^{l_k}\to 
BCT^{l_1 +...+l_k},$$
and
$$d^{k,l_1 ,...,l_k}:BCT^{k}\times BT^{l_1}\times ...\times BT^{l_k}\to 
BCT^{l_1 +...+l_k}.$$

An {\em $A_n$-map} between $A_n$-spaces $X$ and $Y$ is a collection of 
continuous maps $BCT^k\to Hom(X^k ,Y)$ for $k\leq n-1$, which are compatible 
with the $A_n$-structures. 

For example an $A_3$-map between $A_2$-spaces $X,Y$ is a continuous map $f:X\to 
Y$, which is the image of $pt = BCT^1 \to Hom(X,Y)$ and a homotopy $\alpha: 
f\mu\to \mu (f\times f)$, which is defined by the map $I = BCT^2 \to 
Hom(X^2 ,Y)$ ($f\mu$ is the image of the 0-cell which corresponds to the 
$1$-labeled tree $t(2)$, $\mu (f\times f)$ is the image of the 0-cell, 
corresponding to the $0$-labeled tree $t(2)$ and homotopy $\alpha$ is the image 
of the 1-cell, corresponding to the $\epsilon$-labeled tree $t(2)$).
We will call this structure also an $A_3$-structure on the map $f$. Similarly 
$A_n$-structure on the map $f:X\to Y$ between $A_n$-spaces is an $A_n$-map 
structure such that the image of $pt = BCT^1 \to Hom(X,Y)$ coincides with $f$.  

Let us note that $A_3$ map between $A_3$-spaces $X,Y$ allows to define the 
structure of $A_3$-bimodule over an $A_3$-space $X$ on a space $Y$. In 
particular if $Y$ is a loop, space we have defined topological cohomology 
$H^* (X,\Omega^* Y)$.  

\begin{th}\label{map}
Let $(X,\mu_X ),(Y,\mu_Y )$ be $A_n$-spaces (and $Y$ also independently a loop 
spaces if $n=3$) and let $f:X\to Y$ be a $A_n$-map. 
Then the images $f_* (\alpha (X)), f^* (\alpha (Y))$ of their 
canonical classes in $H^{n+1}(X,\Omega^{n-2}Y)$ coincide.  

If the $A_n$-structures on $X,Y$ can be extended to an $A_{n+1}$-structure, then 
the class $\alpha (f)$ in the cohomology $H^{n}(X,\Omega^{n-2}Y)$ is defined, 
does not depend of the choice of $n$-dimesional component of $A_n$-structure on  
$f$ and is trivial iff this structure can be extended to an $A_{n+1}$-structure.
\end{th}
Proof:

Let us note that an $A_n$-map structure on $f$ (and composition operations on 
$BCT^*$ and $Hom(X^* ,Y)$) allows us to define maps $BCT^{l}_{\leq n-2}\to 
Hom(X^l ,Y)$ for any $l$. Here $BCT^{l}_{\leq n-2}$ denotes the union of 
$i$-dimensional cells in $BCT^l$ with $i\leq n-2$. 

Let us consider the element $\alpha (f)$ in $[X^{n},\Omega^{n-2}Y]$ which 
corresponds to the map of $n-2$-sphere $\partial BCT^{n} = 
BCT^{n}_{\leq n-2}\to Hom(X^{n},Y)$. This element satisfies the equation, which 
is defined by the map $BCT^{n+1}_{\leq n-2}\to Hom(X^{n+1},Y)$. The cell complex 
$BCT^{n+1}_{\leq n-2}$ is 
a union of boundaries of $n-1$-cells, which correspond to tricolored trees with 
$n+1$ ends and $n_{\epsilon}+1$ vertices (which means that all exept one 
vertices of the tree are $\epsilon$-labeled). 
Nondegenerated components correspond to the trees with $n_{\epsilon}\leq 1$. 
The case $n_{\epsilon}=0$ consists of two $n-2$-spheres (two tricolored trees of 
the form $t(n+1)$: one is $1$-labeled and another $0$-labeled) which classes 
in $[X^{n+1},\Omega^{n-2}Y]$ are $f_* (\alpha (X))$ and $f^* (\alpha (Y))$ 
respectively. 
In the case $n_{\epsilon}=1$, we have a trees of the form 
$\partial^{i}_{n}(t(n))$ for $i = 0,...,n+1$. Nondegenerate colored trees  
correspond to the case where the tree $t(n)$ is $\epsilon$-labeled. 
Their classes in $[X^{n+2},\Omega^{n}X]$ are equal to 
$\partial^{i}_{n+1}(\alpha (X))$ respectively. Hence the map 
$BT^{n+2}_{\leq n-2}\to Hom(X^{n+2},X)$ defines the equation $\partial 
(\alpha (X)) = 0$ in $[X^{n+2},\Omega^{n-2}X]$. 
$\Box$

Let us note that for an $A_4$-space $(X,\mu_X )$ and $A_4$-bymodule $(Y,\mu_Y )$ 
the collection of 
abelian groups $H^{p}(X,\Omega^{q}Y)$ are equipped with a differential 
$d^{p,q}:H^{p}(X,\Omega^{q}Y)\to H^{p+2}(X,\Omega^{q+1}Y)$. 
Indeed, let $\chi$ be a representative for some cocycle $x\in 
C^{p}(X,\Omega^{q}Y)$. We can choose the homotopy $\tau$ between the map 
$\partial (\chi ):X^{p+1}\to \Omega^{q}Y$ and the trivial map $1_{p+1}$ (the map 
into fixed point). Since $\partial (1_{p+1})$ and $\partial^{2}(\chi )$ can be 
canonically identified with the trivial map $1_{p+2}:X^{p+2}\to\Omega^{q}Y$ the 
homotopy $\partial (\tau )$ can be regarded as an autohomotopy of $1_{p+2}$ or 
as a map $K^{p+2}\to\Omega^{q+1}Y$. It is easy to see that the class of this map 
in $[X^{p+2},\Omega^{q+1}Y]$ is a cocycle and the class of this cocycle in 
$H^{p+2}(X,\Omega^{q+1}Y)$ does not depend of the choice of $\chi$ and $\tau$ 
and defines a homomorphism $d:H^{p}(X,\Omega^{q}Y)\to H^{p+2}(X,\Omega^{q+1}Y)$. 

Direct checking shows that the canonical class of an $A_n$-space ($A_n$-map) 
lies in the cernel of the differential $d$ and the image of the canonical class 
in $E_{3}^{p,q} = ker(d^{p,q})/im(d^{p-2,q-1})$ does not depend on the 
choice of not only the $n$-component of the $A_n$-structure (theorems 
\ref{space},\ref{map}) by also the $n-1$-component. 

Moreover for an $A_n$-space $(X,\mu_X )$ and $A_n$-bymodule $(Y,\mu_Y )$ there 
is defined a "restricted spectral sequence" consisting of abelian groups 
$E_{k}^{p,q}$ ($k\leq n$) and differentials  
$d^{p,q}_{k}:E_{k}^{p,q}\to E_{k}^{p+k,q+k-1}$ ($k\leq n-1$) such that 
$E_{k}^{p,q} = ker(d^{p,q}_{k-1})/im(d^{p-k,q-k+1}_{k-1})$. 
This restricted spectral sequence can be constructed using filtrations on the 
groups $[X^* ,Y]$ which are defined by the following notion. 

An {\em $A_k$-cocycle} of an $A_n$-space $X$ with coefficients in an 
$A_n$-bymodule $Y$ ($k\leq n$) is the class in $[X^m ,Y]$ of a map $f:X^m \to Y$ 
with the homotopies $\tau_1 :\partial (f)\to 1$, 
$\tau_i :\partial (\tau_{i-1})\to 1$ for $2\leq i\leq k$.

It is useful for applications to change the group of values of canonical 
classes. 
Since an $A_3$-space (bimodule) structure on a loop space induces a graded ring 
(bimodule) structure on its homotopy groups, the natural map from 
$C_{*}(X,Y)$ to the Hochschild complex 
$C_{*}(\pi_{\cdot}(X),\pi_{\cdot}(Y))$ of the ring $\pi_{*}(X)$ with 
coefficients in the bimodule $\pi_{*}(Y)$ induces the homomorphism of 
cohomology
$$H^{*}(X,Y) \longrightarrow HH_{*}(\pi_{*}(X),\pi_{*}(Y)).$$
So we have the class $\alpha (X)\in HH_{n+3}(\pi_{0}(X),\pi_{n}(X))$ for any 
$A_n$-space $X$. 

\section{Cohomological obstruction for monoidality}\label{com}
In \cite{qui} Quillen associated to an abelian (exact) category ${\cal A}$ a 
topological space $BQ{\cal A}$ such that exact functors between categories 
define continuous maps between the corresponding spaces and isomorphisms of 
functors define homotopies between the corresponding maps. In other words, 
Quillen's space is a 2-functor from the 2-category of abelian (exact) categories 
(with isomorphisms of functors as 2-morphisms) to the 2-category of topological 
spaces.
Waldhausen \cite{wld} proved that the 2-functor $K = \Omega BQ$ is permutative 
(in some sense) with the natural product operations. Namely, he constructed a  
continuous map 
$K({\cal A})\wedge K({\cal B})\rightarrow K({\cal C})$ for any biexact functor 
${\cal A}\times {\cal B}\rightarrow {\cal C}$. He also proved that $K({\cal A})$ 
is an infinite loop space for any abelian category $\cal A$. 
The homotopy groups $K_{*}({\cal A}) = \pi_{*}({\cal A})$ of the Waldhausen 
space $K({\cal A})$ are called the {\em algebraic K-theory} of the category 
${\cal A}$.
The 2-categorical nature of the functor $K$ implies that the Waldhausen space 
$K({\cal G})$ is an $A_{\infty}$-space for any abelian monoidal category 
$\cal G$
with respect to the product corresponding to the monoidal product in $\cal G$. 
In particular, the space $K({\cal M}_{k}(S))$ is $A_{\infty}$-space for any 
monoidal (over the field $k$) fusion algebra $S$. 

The structure of the space $K({\cal M}_{k}(S))$ motivates the following 
construction. 
If $S$ is a fusion ring and $K$ is a homotopy associative space (which is also 
a loop space), then 
we can construct a new space $K(S) = \bigvee_{s\in S}K$ with the product 
defined by the fusion rule of $S$. Namely the component of $\mu_{K(S)}$ 
between the $i$ and $j$-labeled components of $K(S)$ to the $k$-labeled 
has the form ${\mu_{K(S)}}_{i,j}^{k} = m_{i,j}^{k}\mu_{K}$ where the 
multiplication by the non-negative integer $m_{i,j}^{k}$ is defined by means of 
internal loop-space structure. 

We say that the fusion ring $S$ is of {\em $A_n$-type} with respect to an 
$A_{\infty}$-space $K$ if the natural $A_3$-structure on $K(S)$ is extendeble 
to an $A_n$-structure. 
The morphism of fusion rings $S,S'$ is of {\em $A_n$-type} with respect to 
$A_{\infty}$-space $K$ if the corresponding map $K(S)\to K(S')$ has an 
$A_n$-structure. 

The {\em canonical} class $\alpha (S)$ of the fusion ring $S$ with respect to 
the $A_{\infty}$-space $K$ is a class of the homotopy associative space $K(S)$ 
in the 
Hochschild cohomology $HH^* (A(S),A(S)\otimes\pi_* (K))$. The next theorem is a 
direct corollary of the definitions.
\begin{th}
For the monoidal (over the field $k$) fusion ring $S$, the canonical cohomology 
class $\alpha_* (S)\in HH^* (A(S),A(S)\otimes K_* (k))$ is trivial. 
\end{th}

Properties of canonical classes of $A_n$-spaces proved in the previous section 
(theorem \ref{map}) imply the following theorem. 
\begin{th}
Let the fusion rings $S,S'$ be of $A_n$-type with respect to $A_{\infty}$-space 
$K$. 
Let $f:A(S)\to A(S')$ be the homomorphism of enveloping rings which corresponds 
to the morphism of fusion rings $S$ and $S'$ of $A_n$-type. Then the canonical 
classes of these fusion rings satisfy the following condition:
$$f_* (\alpha (X)) = f^* (\alpha (Y)) 
\in HH^{n+3}(A(S),A(S')\otimes\pi_{n}(K)).$$ 
If the fusion rings $S,S'$ are of $A_{n+1}$-type then the class $\alpha (f)$ in 
cohomology $HH^{n+2}(A(S),A(S')\otimes\pi_{n}(K))$ is defined and is trivial if 
the map $f$ admits a $A_{n+1}$-structure.
\end{th}
For example any fusion ring is of $A_3$-type with respect to arbitrary an 
$A_{\infty}$-space $K$. Hence the class $\alpha (S)\in 
HH^{4}(A(S),A(S')\otimes\pi_{1}(K))$ is defined. Any morphism of fusion rings 
is also of $A_3$-type, hence the images $f_* (\alpha (X)), f^* (\alpha (Y))$ in 
$HH^{4}(A(S),A(S')\otimes\pi_{1}(K))$ of the first classes of the fusion rings 
$S,S'$, connected by the morphism $f$, coincide. 

The first class with respect to the space $K({\cal M}(k)) = BGL(k)^+$ for the 
field $k$ admits a more direct describtion. Namely it coincides with the class 
of the Hochcshild cocycle $A\in Z^4 (A(S),A(S)\otimes k^* )$, which is defined 
by the following
$$A(x_1 ,x_2 ,x_3 ,x_4 ) = \sum_{x\in S}\det (A^{x}_{x_1 ,x_2 ,x_3 ,x_4}),$$ 
where $k^*$ is the group of invertible elements of the field $k$ and 
$det (A^{x}_{x_1 ,x_2 ,x_3 ,x_4})$ is the determinant of the linear automorphism, which is a 
(clockwise) circumference of the diagram $A^{x}_{x_1 ,x_2 ,x_3 ,x_4})$. 

It follows from the describtion that first class of fusion ring $S$ in 
$HH^{4}(A(S),A(S)\otimes k^* )$ is the image of the class in 
$HH^{4}(A(S),A(S)\otimes\pi_1 (B\Sigma_{\infty}^{+}))$ where the space 
$B\Sigma_{\infty}^{+}$ is Quillen's plus-construction applyed to the 
infinite symmetric group and (an analog of the Waldhausen space $K({\cal A})$ 
for the category of finite sets). Barrat-Priddy-Quillen 
theorem \cite{bp,seg1,seg2} states that this space is homotopy equivalent to the 
spheric spectrum. In particular its homotopy groups (stable homotopy groups of 
spheres) are torsion and 
$$\pi_1 (B\Sigma_{\infty}^{+}) = {\bf Z}/2{\bf Z},\quad 
\pi_2 (B\Sigma_{\infty}^{+}) = {\bf Z}/2{\bf Z},\quad 
\pi_3 (B\Sigma_{\infty}^{+}) = {\bf Z}/24{\bf Z}.$$

\section{Calculation of the first obstruction}\label{calc}
In this section we will give the combinatorial describtion of the cocycle 
representing the first obstruction and calculate its class in the case of two 
element fusion ring with identity. 

To write the first obstruction in explicit form, choose a linear order $>$ on 
the set $S$ and linearly ordered sets 
$$X^{x}_{y,z} = \{ f^{x}_{y,z}(i),\quad i = 1,...,m^{x}_{y,z}\}$$ 
for any $x,y,z\in S$. This data allows us to define well order on the sets 
$$\bigsqcup_{u\in S}X^{x}_{y,u}\times X^{u}_{z,w},\quad  
\bigsqcup_{v\in S}X^{v}_{y,z}\times X^{x}_{v,w}$$  
and to define the map $\Phi^{x}_{y,z,w}$ between them as a (unique) order 
preserving bijection. More precisely elements of the first set
$$f^{x}_{y,z,w}(u,i,j) = f^{x}_{y,u}(i)\times f^{u}_{z,w}(j)$$
parametrize by the collections 
$$(u,i,j),\quad u\in S, i=1,...,m^{x}_{y,u}, j=1,...,m^{u}_{z,w},$$ 
which are ordered lexicographically
$$(u,i,j)>(u',i',j')\Longleftrightarrow \left\{
\begin{array}{ccc}
u>u' \\ u=u', & i>i' \\ u=u', & i=i', & j>j' 
\end{array}\right. .$$
The elements of the second set
$$g^{x}_{y,z,w}(v,s,t) = f^{v}_{y,z}(t)\times f^{x}_{v,w}(s)$$
which are parametrized by the collections 
$$(v,s,t),\quad v\in S, s=1,...,m^{v}_{y,z}, j=1,...,m^{x}_{v,w},$$ 
ordered analogously:  
$$(v,s,t)>(v',s',t')\Longleftrightarrow \left\{
\begin{array}{ccc}
v>v' \\ v=v', & s>s' \\ v=v', & s=s', & t>t' 
\end{array}\right. .$$
The following lemma will be useful for the calculation of the first obstruction.
\begin{lem}\label{sign}
1. The sign of unique order preserving bijection between lexicographically 
ordered sets $X\times Y$ and $Y\times X$ (the permutation of factors) equals 
$(-1)^{{|X|\choose 2}{|Y|\choose 2}}$.

2. The sign of unique order preserving bijection between lexicographically 
ordered disjoint unions $\bigsqcup_{a,b\in S}X(a,b)$ and 
$\bigsqcup_{b,a\in S}X(a,b)$ equals 
$$(-1)^{\sum_{a_1 > a_2 ,b_1 < b_2}|X(a_1 ,b_1 )||X(a_2 ,b_2 )|}.$$ 
\end{lem}
Proof:

The sing the automorphism of the set $X$ sending one well order $>_1$ to another 
$>_2$ is determined by the evenness of the number of inversions of orders, e.g. 
the quantaty of pairs $x,y\in X$ such that $x>_1 y$ and $x<_2 y$. 

1. The inversion of two lexicographic orders on $X\times Y$ is a pair 
$(x_1 ,y_1 ),(x_2 ,y_2 )\in X\times Y$ such that $x_1 >x_2$ and $y_1 <y_2$. 
The quantaty of these pairs is ${{|X|\choose 2}{|Y|\choose 2}}$.

2. The inversion of two lexicographic orders on disjoint union 
$\bigsqcup_{a,b\in S}X(a,b)$ is a pair $x\in X(a_1 ,b_1 ), y\in X(a_2 ,b_2 )$ 
such that $a_1 >a_2$ and $b_1 <b_2$. The quantaty of these pairs is 
${\sum_{a_1 > a_2 ,b_1 < b_2}|X(a_1 ,b_1 )||X(a_2 ,b_2 )|}$. $\Box$

To write down the first obstruction let us fix an order on the set of vertices 
of any planar binary tree. 

First of all we can regard any planar tree as a partial order on the set of its 
vertices $V$. Thus the set $V$ decomposes into the disjoint union of subsets of 
uncomparible vertices. Since the tree is planar we can draw all vertices from 
these subsets lying on the same horizontal line. We can define the well order 
on $V$ ordering them it the left to right. 

The order on the set of vertices of the planar binary tree $T$ with $n$ 
ends allows to define an (lexicographic) order on the set 
$$X^{y}_{x_1 ,...,x_n}(T) = 
\bigsqcup_{f:E\to S}\times_{v\in V}X^{f(e_1 )}_{f(e_2 ),f(e_3 )},$$
where sum is over all marking of the tree $T$, e.g. the functions from the set 
of edges $E$ to the fusion ring $S$, such that the values on the ends are 
$x_1 ,...,x_n$; $e_1\in E$ is unique edge with the end $v$ and $e_2 ,e_3$ are 
edges with the begining $v$. 

In particular we can define the orders on the (five) sets corresponding to the 
planar binary trees with four ends which parametrize the vertices of the 
pentagon diagram $A^{x}_{x_1 ,x_2 ,x_3 ,x_4}$:
$$\bigsqcup_{a,b\in S}X^{x}_{x_1 ,a}\times X^{a}_{x_2 ,b}\times X^{b}_{x_3 
,x_4},\quad 
\bigsqcup_{c,b\in S}X^{x}_{c,b}\times X^{c}_{x_2 ,x_2}\times X^{b}_{x_3 ,x_4},$$
$$\bigsqcup_{a,e\in S}X^{x}_{x_1 ,a}\times X^{a}_{e,x_4}\times X^{e}_{x_2 ,x_3}, 
\bigsqcup_{d,e\in S}X^{x}_{d,x_4}\times X^{d}_{x_1 ,e}\times X^{e}_{x_2 ,x_3},$$  
$$\bigsqcup_{c,d\in S}X^{x}_{d,x_4}\times X^{d}_{c,x_3}\times X^{c}_{x_1 
,x_2}.$$ 
These orders allow to correspond the (order preserving) map to any arrow of the 
diagram $A^{x}_{x_1 ,x_2 ,x_3 ,x_4}$. By the other hand we can define a map 
corresponding to any arrow of the diagram $A^{x}_{x_1 ,x_2 ,x_3 ,x_4}$ by means 
of the maps $\Phi$. These pairs of maps does not coincide with each other. Using 
lemma \ref{sign} we can calculate the signs of the differences: 
$$S^{x}_{x_1 ,x_2 ,x_3 ,x_4} = $$
$$\sum_{a>a',b<b'}m^{x}_{x_1 ,a}m^{a}_{x_2 ,b}m^{b}_{x_3 ,x_4}
m^{x}_{x_1 ,a'}m^{a'}_{x_2 ,b'}m^{b'}_{x_3 ,x_4} +$$
$$\sum_{c>c',b<b'}m^{x}_{c,b}m^{c}_{x_1 ,x_2}m^{b}_{x_3 ,x_4}
m^{x}_{c',b'}m^{c'}_{x_1 ,x_2}m^{b'}_{x_3 ,x_4} +$$
$$\sum_{a>a',e<e'}m^{x}_{x_1 ,a}m^{a}_{e,x_4}m^{e}_{x_2 ,x_3}
m^{x}_{x_1 ,a'}m^{a'}_{e',x_4}m^{e'}_{x_2 ,x_3} +$$
$$\sum_{d>d',e<e'}m^{x}_{d,x_4}m^{d}_{x_1 ,e}m^{e}_{x_2 ,x_3}
m^{x}_{d',x_4}m^{d'}_{x_1 ,e'}m^{e'}_{x_2 ,x_3} +$$
$$\sum_{d>d',c<c'}m^{x}_{d,x_4}m^{d}_{c, x_3}m^{c}_{x_1 ,x_2}
m^{x}_{d',x_4}m^{d'}_{c', x_3}m^{c'}_{x_1 ,x_2} +$$
$$\sum_{c,b}m^{x}_{c,b}{m^{c}_{x_1 ,x_2}\choose 2}{m^{b}_{x_3 ,x_4}\choose 2}.$$
The cocycle representing the class of the first obstruction can be written as 
follows
$$\alpha (x_1 ,x_2 ,x_3 ,x_4) = 
\sum_{x\in S}x\otimes\alpha^{x}_{x_1 ,x_2 ,x_3 ,x_4},$$
where $\alpha^{x}_{x_1 ,x_2 ,x_3 ,x_4} = (-1)^{S^{x}_{x_1 ,x_2 ,x_3 ,x_4}}$. 
\bigskip

Example. 

Let us consider the case of two elements fusion ring $S = \{ e,x\}$ with the 
fusion rule $x*x = mx + n$. Let us note that the first obstruction lies in 
$H^4 (A(S),A(S)\otimes {\bf Z}/2{\bf Z})\cong H^4 (A(S),A(S)/2A(S))$. 
\begin{lem}
The fourth Hochcshild cohomology group $H^4 (A(S),A(S)/2A(S))$ 
of the above algebra $A(S)$ are isomorphic to 
$$A(S)/(2,m)A(S)\cong\left\{
\begin{array}{cc}
0, & m\equiv 1(2) \\
{\bf Z}/2{\bf Z}\oplus{\bf Z}/2{\bf Z}, & m\equiv 0(2) 
\end{array}\right. .$$
The isomorphism sends the class of the cocyle $\alpha$ of the standart complex 
to the class of its value $\alpha (x,x,x,x)\in A(S)/(2,m)A(S)$. 
\end{lem}
Proof:
Let us use the presentation of Hochcshild cohomology groups $H^k (A(S),M)$ with 
coefficients in the $A(S)$-bimodule $M$ as the groups of extensions 
$Ext^{k}_{A(S)-A(S)}(A(S),M)$ in the category of $A(S)$-bimodules. 
Since $A(S) = {\bf Z}[x,x^2 =mx+n]$ we can identify free $A(S)$-bimodule of 
rank one $A(S)\otimes_{\bf Z}A(S)$ with ${\bf Z}[x_1 ,x_2 , x_{i}^2 = mx_i +n]$.
The $A(S)$-bimodule $A(S)$ has a 2-periodic resolution 
$$A(S)\gets^{\mu}A(S)\otimes A(S)\gets^{x_1 -x_2}A(S)\otimes A(S) 
\gets^{x_1 +x_2 -m}$$
$$A(S)\otimes A(S)\gets^{x_1 -x_2}A(S)\otimes A(S)\gets^{x_1 +x_2 -m}...$$
Applying the functor $Hom_{A(S)-A(S)}(?,M)$ to this resolution we can 
calculate the cohomology $H^* (A(S),M)$  
$$H^{2i}(A(S),M)\equiv\{y\in M, xy=yx\}/\{xz+zx-mz, z\in M\},$$
$$H^{2i+1}(A(S),M)\equiv\{y\in M, xy+yx=my\}/\{xz-zx, z\in M\}.$$
Since in our case bimodule $M = A(S)/2A(S)$ is symmetric we have 
$$H^4 (A(S),A(S)/2A(S))\equiv A(S)/(2,m)A(S).$$
$\Box$

The unique nontrivial factor of the expression $\alpha^{e}_{x,x,x,x}$ is the 
last one which equals $(-1)^{{n\choose 2} + n{m\choose 2}}$. 
The last factor of the expression $\alpha^{x}_{x,x,x,x}$ is equal to 
$(-1)^{m{m\choose 2}}$ and all other coincide with $(-1)^{nm}$. Hence 
$$\alpha ^{x}_{x,x,x,x} = (-1)^{m{m\choose 2} + nm}.$$

It follows from the above describtion that the first obstraction of the fusion 
ring $S$ is nontrivial iff 
$$m\equiv 0(2), n\equiv 2,3(4),\quad\mbox{or}\quad m\equiv 2(4), n\equiv 1(4).$$

\end{document}